\documentclass[12pt]{amsart}
\usepackage{amssymb}
\usepackage[margin = 1in]{geometry}
\usepackage[shortlabels]{enumitem}
\newcommand{\norm}[1]{\lVert#1\rVert}
\newtheorem{theorem}{Theorem}
\newtheorem{lemma}{Lemma}
\theoremstyle{definition}

\newtheorem{remark}{Remark}
\newtheorem{example}{Example}
\newtheorem{assumption}{Assumption}
\begin{document}
\author{Andrey Sarantsev}
\title[Convergence Rate in Wasserstein Distance]{Convergence rate to equilibrium in Wasserstein distance for reflected jump-diffusions}
\address{University of Nevada in Reno, Department of Mathematics and Statistics}
\email{asarantsev@unr.edu}

\begin{abstract}
Convergence rate to the stationary distribution for continuous-time Markov processes can be studied using Lyapunov functions. Recent work by the author provided explicit rates of convergence in special case of a reflected jump-diffusion on a half-line. These results are proved for total variation distance and its generalizations: measure distances defined by test functions regardless of their continuity. Here we prove similar results for Wasserstein distance, convergence in which is related to convergence for continuou test functions. In some cases, including the reflected Ornstein-Uhlenbeck process, we get faster exponential convergence rates for Wasserstein distance than for total variation distance. 
\end{abstract}
\maketitle
\thispagestyle{empty}

\section{Introduction} Consider a Markov process $X = (X(t),\, t \ge 0)$ on a metric state space $S$ with transition function $P^t(x, \cdot)$. Assume it has a unique stationary distribution, or invariant measure $\pi$: If $X(0) \sim \pi$ then $X(t) \sim \pi$ for all $t \ge 0$. We are interested in convergence 
\begin{equation}
\label{eq:conv}
P^t(x, \cdot) \to \pi(\cdot),\quad t \to \infty.
\end{equation}
In which norm does it hold? How fast is this convergence? One well-known tool is a Lyapunov function. Define the $g$-norm for signed measures on $S$ as follows: 
\begin{equation}
\label{eq:g-norm}
\norm{\nu}_g := \sup\limits_{|f| \le g}|(\nu, f)|,\quad (\nu, f) := \int_Sf(x)\nu(\mathrm{d}x).
\end{equation}
For $g = 1$, this becomes a total variation norm $\norm{\cdot}_{\mathrm{TV}}$. Assume $\mathcal L$ is the generator of this process. If there is a function $V : S \to [1, \infty)$ such that for a compact set $K \subseteq S$,
\begin{equation}
\label{eq:LF}
\mathcal L V(x) \le -kV(x),\quad x \in S\setminus K,
\end{equation}
called a {\it Lyapunov function}, then (under additional technical assumptions) the convergence in~\eqref{eq:conv} is exponential:  
\begin{equation}
\label{eq:conv-exp}
\norm{P^t(x, \cdot) - \pi(\cdot)}_{\mathrm{TV}} \le C(x)e^{-\varkappa t}
\end{equation}
for some constant $C(x)$ depending on $x \in S$, and $\varkappa > 0$. This statement can also be proved for the $V$-norm $\norm{\cdot}_V$ from~\eqref{eq:g-norm}, for the same Lyapunov function $V$, see \cite{BCG08, DMT95, MT93a, MT93b}. However, to find or estimate an exact constant $\varkappa$ is challenging: \cite{MT94, MCMC}. But we can conclude that
$\varkappa = k$ for  $S = [0, \infty),\, K = \{0\}$, and the process $X$ is stochastically ordered, \cite{LMT96, Sar16}. These results are applied to  reflected diffusions and jump-diffusions. A reflected diffusion can be informally described informally as follows: Inside the half-line, it behaves between jumps as a solution to one-dimensional stochastic differential equation with drift $g$ and diffusion $\sigma^2$:
$$
\mathrm{d}X(t) = g(X(t))\,\mathrm{d} t + \sigma(X(t))\,\mathrm{d}W(t),
$$
where $W$ is a Brownian motion. When it hits $0$, it is reflected back to the half-line. A reflected jump-diffusion behaves as a reflected diffusion but, in addition, it can jump with a certain rate. The jump destination is also random, and there are finitely many jumps in finite time. We applied our results to queueing theory \cite{Q} and risk theory \cite{Risk}. 

For the proof, we use the coupling method: Take two copies $X_1$ and $X_2$ of this process starting from points $X_2(0) \ge X_1(0) \ge 0$. By stochastic ordering, we can couple them (that is, create them on the same probability space) so that $0 \le X_1(t) \le X_2(t)$. Take $\tau$ to be the hitting moment of zero by $X_2$. Then by stochastic ordering $X_1(\tau) = 0$, and we can assume $X_1(t) = X_2(t)$ for $t > \tau$. This enables us to estimate the distance between $P^t(x_1, \cdot)$ and $P^t(x_2, \cdot)$. This coupling method is commonly used for convergence proofs, see articles \cite{Hairer, Lindvall83, Lindvall86} and the book \cite{LindvallBook}. In our case, when $S = \mathbb R_+$ and $K = \{0\}$, this method is particularly powerful.

In this short note, we apply this method to Wasserstein distance instead of $\norm{\cdot}_V$-norm. Wasserstein distance is defined via optimal couplings of measures; one can think of it as ``earth mover'' distance: If we have two piles of sand with equal volume, how much work does it require to move the first pile in the place of the second pile? This distance is related to {\it optimal transport} problems, see the fundamental monograph \cite{Villani}, in particular Chapter 6. As noted there, convergence in Wasserstein distance of order $p \ge 1$ is equivalent to uniform boundedness of $p$th moments and weak convergence of measures. The latter, in turn, is equivalent to convergence for continuous bounded test functions. Thus Wasserstein distance is fundamentally different from total variation norm or other norms as in~\eqref{eq:g-norm}, which use all test functions $f : S \to \mathbb R$, not just continuous ones. Previous work on convergence in this distance is scant, see \cite{Oleg}. In this article, we find convergence rates and compare them with that for total variation and other similar distances. For some processes, including a reflected Ornstein-Uhlenbeck process, convergence is faster in Wasserstein distance than in total variation or the $\norm{\cdot}_V$-norm.

We note that the interplay between Wasserstein distance and Lyapunov functions is remarkable. Concentration of measure Talagrand inequalities for a distribution $\mathbb P$ compare Wasserstein distance from $\mathbb P$ to a test measure 
$\mathbb Q$ with relative entropy of $\mathbb Q$ with respect to $\mathbb P$. These inequalities were found for common stochastic processes on $[0, T]$; see \cite{Pal, CG, PS19, Samson} and references therein. These concentration inequalities are related to other functional inequalities for finite time horizon, which can be obtained using Lyapunov functions, \cite{BCG08}. But these articles do not address the long-term behavior of stochastic processes, even if these processes are convergent in the long run (like an Ornstein-Uhlenbeck process). In this work, we fill this gap and combine Lyapunov functions with Wasserstein distance framework to study rates of long-term convergence.

This short note is organized as follows. In Section 2, we define all notation and concepts. In Section 3, we state our main results and give examples. Section 4 is devoted to proofs. 

\section{Notation, Definitions, and Background}

Define the {\it Skorohod space} $\mathbb D([s, t])$ of right-continuous function with left limits $[s, t] \to \mathbb R$, with the distance $\rho_{\mathbb D}$ satisfying the following estimate:
\begin{equation}
\label{eq:Skorohod}
\rho_{\mathbb D}(x, y) \le \sup\limits_{s \le u \le t}|x(u) - y(u)|,\quad x, y \in D[s, t].
\end{equation}
For two probability measures $\mathbb P$ and $\mathbb Q$ on the same metric space $(S, \rho)$, and for a  $p \ge 1$, {\it Wasserstein distance of order} $p$ is defined as
$$
\mathcal W_p(\mathbb P, \mathbb Q) = \inf\left[\mathbb E \, \rho(X, Y)^p \right]^{1/p},
$$
where the infimum is taken over all couplings $(X, Y)$ such that $X \sim \mathbb P$ and $Y \sim \mathbb Q$. A family of finite measures $(Q_x)_{x \ge 0}$ on $\mathbb R_+$ is {\it stochastically ordered} if
$$
Q_x(\mathbb R_+) = \mathrm{const},\quad Q_{x}([z, \infty)) \le Q_y([z, \infty)),\quad 0 \le x \le y,\quad z \ge 0.
$$
We operate on a filtered probability space $(\Omega, \mathcal F, (\mathcal F_t)_{t \ge 0}, \mathbb P)$. Consider a reflected jump-diffusion on $\mathbb R_+$ with drift $g : \mathbb R_+ \to \mathbb R$, diffusion $\sigma : \mathbb R_+ \to \mathbb R_+$, and family of jump measures $(\nu_x)_{x \ge 0}$ on $\mathbb R_+$. Take an $(\mathcal F_t)_{t \ge 0}$-Brownian motion $W = (W(t),\, t \ge 0)$. 

We formally define a reflected diffusion without jumps, formalizing the informal description from the Introduction. This is an a.s. continuous $(\mathcal F_t)_{t \ge 0}$-adapted $\mathbb R_+$-valued process $X = (X(t),\, t \ge 0)$ such that there exists an adapted continuous nondecreasing process $\ell = (\ell(t),\, t \ge 0)$ with $\ell(0) = 0$, which can increase only when $X(t) = 0$, such that 
\begin{equation}
\label{eq:SDER}
X(t) = X(0) + \int_0^tg(X(s))\,\mathrm{d}s + \int_0^t\sigma(X(s))\,\mathrm{d}W(s) + \ell(t).
\end{equation}
Add jumps to~\eqref{eq:SDER} to get a reflected jump-diffusion by {\it piecing out:} At a point $x \in \mathbb R_+$ it jumps with rate $N(x) := \nu_x(\mathbb R_+)$, with jump destination distributed as $N^{-1}(x)\nu_x(\cdot)$. This is done by {\it piecing out:} Assuming that $N(x) = \Lambda$ is constant (as in Assumption 1 below),  we fix jump times $\tau_1 < \tau_2 < \ldots$ as Poisson point process with intensity $\Lambda$: That is, $\tau_{k+1} - \tau_k \sim \mathrm{Exp}(\Lambda)$ i.i.d. with convention $\tau_0 := 0$. Finally, we run this process as a reflected diffusion without jumps between $\tau_k$ and $\tau_{k+1}$. More details are in \cite[Section 2]{Sar16} and in \cite{Piece}.

\begin{assumption} The functions $g$ and $\sigma$ are continuous, this family of measues is weakly continuous: $\nu_y \to \nu_x$ as $y \to x$ weakly, and $N(x) = \Lambda$ is a constant.
\end{assumption}

\begin{remark} It suffices to assume $\sup_{x \ge 0}N(x) = \overline{N} < \infty$ instead of constancy from Assumption 1. Indeed, we can replace measures $\nu_x$ with $\nu'_x := \nu_x + (\overline{N} - N(x))\delta_x$. Indeed, jumping from $x$ to $x$ means not jumping at all, and new measures $\nu'_x$ satisfy $\nu'_x(\mathbb R_+) = \overline{N} = \mathrm{const}$.
\end{remark} 

It was shown in \cite{Sar16} that under Assumption 1, there exists in the weak sense a unique in law version of this process for any initial condition $X(0)$ (which can be random). This process is Feller continuous strong Markov, with generator $\mathcal Lf$ for $f \in C^2(\mathbb R_+)$ with $f'(0) = 0$:
$$
\mathcal Lf(x) = g(x)f'(x)+ \frac{\sigma^2(x)}2f''(x) + \int_0^{\infty}[f(y) - f(x)]\,\nu_x(\mathrm{d}y),
$$
If the initial distribution is $X(0) \sim \rho$, then $X(t) \sim \rho P^t$ with $\rho P^t$ for $t \ge 0$ defined as
$$
\rho P^t(B) := (\rho, P^t(\cdot, B)) = \int_0^{\infty}\rho(\mathrm{d}x)\,P^t(x, B),\quad B \subseteq [0, \infty).
$$
For example, if $x \ge 0$, then $\delta_x P^t(B) = P^t(x, B)$. A probability distribution $\pi$ on $\mathbb R_+$ is called a {\it stationary distribution} for this Markov process if $\pi = \pi P^t$ for $t \ge 0$; or, equivalently, if $X(0) \sim \pi$ implies $X(t) \sim \pi$ for all $t \ge 0$. This Markov process is called {\it stochastically ordered} if for every $t \ge 0$, the family of measures $(P^t(x, \cdot))_{x \ge 0}$ is stochastically ordered. For $0 \le s \le t$ and $x \ge 0$, define $\mathbf{P}^{[s, t]}(x, \cdot)$ to be the distribution in the Skorohod space $D[s, t]$ of the process $X = (X(u),\, u \in [s, t])$ starting from $X(0) = x$.

\section{Main Results}

For a $\lambda > 0$, define $k(\lambda) := -\sup\limits_{x > 0}c(x, \lambda)$, where
\begin{align}
\label{eq:c-k}
\begin{split}
c(x, \lambda) &:= \lambda g(x) + \frac{\lambda^2}2\sigma^2(x) + \int_0^{\infty}\left[e^{\lambda(y-x)} - 1\right]\,\nu_x(\mathrm{d}y),\quad x \in \mathbb R_+.
\end{split}
\end{align}

\begin{assumption} The family of jump measures $(\nu_x)_{x \ge 0}$ is stochastically ordered. 
\end{assumption}

\begin{assumption} The integral in~\eqref{eq:c-k} is well-defined for all $x, \lambda$, and $k(\lambda) > 0$ for some $\lambda > 0$.
\end{assumption}

It was shown in \cite{Sar16} that Assumption 3 holds under the following conditions:
\begin{equation}
\label{eq:mean}
m(x) := g(x) + \int_0^{\infty}(y-x)\,\nu_x(\mathrm{d}y) \le -K_1 < 0,\quad \sigma(x) \le K_2 < \infty,
\end{equation}
for constants $K_1, K_2 > 0$. This can be viewed as ``effective drift'' at point $x \in \mathbb R_+$: The sum of ``true drift'' $g(x)$ and ``implied drift'' from jumps. 

\begin{example} A simple example is a reflected spectrally positive L\'evy process $X$, where $g$ and 
$\sigma$ are constant, and the jump measures $\nu_x$ are defined as follows, for some finite measure $\mu$ on $\mathbb R_+$: $\nu_x([0, x]) = 0$, and $\nu_x([x + z, \infty)) = \mu([z, \infty))$ for all $x, z \ge 0$. This process behaves as a L\'evy process which behaves as a Brownian motion with drift $g$ and diffusion $\sigma^2$ between jumps; jump times behave as a Poisson process with intensity $\mu(\mathbb R_+)$ and jump displacement (to the right) is distributed as the normalized measure $\mu$. As long as it hits zero, it reflects back to the positive half-line. Then the condition~\eqref{eq:mean} becomes
\begin{equation}
\label{eq:meanconst}
g + \int_0^{\infty}z\,\mu(\mathrm{d}z) < 0.
\end{equation}
\end{example}

\begin{example}
Try $g = -1$ and $\sigma = 1$, $\mu \sim \mathrm{Exp}(2)$. Then condition~\eqref{eq:meanconst} holds, and 
$$
-k(\lambda) = c(x, \lambda) = -\lambda + \frac{\lambda^2}2 + \int_0^{\infty}(e^{\lambda z} - 1)\,\mu(\mathrm{d}z) = -\lambda + \frac{\lambda^2}2 + \frac{2}{2 - \lambda} - 1. 
$$
Thus $k$ assumes maximum $k(\lambda_0) = 0.0785$ for $\lambda_0 = 0.304$. 
\end{example}

It was shown (see \cite{Sar16} and references therein) that under Assumptions 1, 2, 3 there exists a unqiue stationary distribution $\pi$ and $(\pi, V) < \infty$ for $V(x) = e^{\lambda x}$; and the $V$-distance between $P^t(x, \cdot)$ and $\pi(\cdot)$ converges to $0$ at least as $e^{-k(\lambda)t}$. We show a similar result for the Wasserstein distance. 

\begin{theorem} Under Assumptions 1, 2, 3, for every $p \ge 1$, for $C := ep/\lambda$, we have
\label{thm:1}
\begin{enumerate}[(a)]

\item For all $x_1, x_2 \ge 0$ and $t \ge 0$, we have:
$$
\mathcal W_p(P^t(x_1, \cdot), P^t(x_2, \cdot)) \le C\cdot\exp\left[p^{-1}\lambda(x_1\vee x_2)\right]\exp\left[- p^{-1}k(\lambda)pt\right]. 
$$
\item  The same estimate holds for path distributions: For $x_1, x_2 \ge 0$ and  $T > t > 0$,
$$
\mathcal W_p\left(\mathbf{P}^{[t, T]}(x_1, \cdot), \mathbf{P}^{[t, T]}(x_2, \cdot)\right) \le  C\cdot\exp\left[p^{-1}\lambda(x_1\vee x_2)\right]\exp\left[-p^{-1}k(\lambda)t\right]. 
$$
\item For probability measures $\rho_1$ and $\rho_2$ on $\mathbb R_+$ and for $t \ge 0$:
$$
\mathcal W_p(\rho_1P^t, \rho_2P^t) \le C\cdot\left[(\rho_1, V) + (\rho_2, V)\right]^{1/p}\exp\left[-p^{-1}k(\lambda)t\right],\quad V(x) := e^{\lambda x}.
$$
\item For every $x \ge 0$ and $t \ge 0$, we get:
$$
\mathcal W_p(P^t(x, \cdot), \pi) \le C\cdot\left[(\pi, V) + V(x)\right]^{1/p}\exp\left[-p^{-1}k(\lambda)t\right]
$$
\end{enumerate}
\end{theorem}

\begin{assumption} The function $\sigma$ is constant. For every $z \in \mathbb R$, the function $x \mapsto \nu_x((x+z, \infty))$ is nonincreasing. There exists a constant $G$ such that $x \mapsto g(x) - Gx$ is nonincreasing.
\end{assumption}

\begin{theorem}
Under Assumptions 1, 2, 3, 4, if $k(\lambda) > pG$, $p \ge 1$, define $K := p^{-1}k(\lambda) - G$. 
\begin{enumerate}[(a)]

\item For all $x_1, x_2 \ge 0,\, t \ge 0$, we have:
$$
\mathcal W_p(P^t(x_1, \cdot), P^t(x_2, \cdot)) \le \exp(p^{-1}\lambda (x_1\vee x_2))|x_1 - x_2|\exp(-Kt).
$$
\item If $G \le 0$, the same estimate holds for path distributions:
$$
\mathcal W_p\left(\mathbf{P}^{[t, T]}(x_1, \cdot), \mathbf{P}^{[t, T]}(x_2, \cdot)\right) \le \exp(p^{-1}\lambda (x_1\vee x_2))|x_1 - x_2|\exp(-Kt).
$$
\item For probability measures $\rho_1, \rho_2$ on $\mathbb R_+$, there exists a constant $C_0 > 0$ such that 
$$
\mathcal W_p(\rho_1P^t, \rho_2P^t) \le C_0\exp(-Kt),\quad V(x) := e^{\lambda x},\quad t \ge 0.
$$
\item For every $t, x \ge 0$, we get:
$$
\mathcal W_p(P^t(x, \cdot), \pi) \le \left[(\pi, V) + V(x)\right]^{1/p}\exp(-Kt).
$$
\end{enumerate}
\label{thm:2}
\end{theorem}

As follows from \cite[Chapter 6]{Villani}, convergence in Wasserstein distance of order $p$ implies convergence of moments up to $p$th order: For every continuous function $f : \mathbb R_+ \to \mathbb R$ with $|f(x)| \le 1 + x^p$ for $x \ge 0$, there exists a constant $C(f) > 0$ such that 
$$
\left|(\nu_1, f) - (\nu_2, f)\right| \le C(f)\mathcal W_p^p(\nu_1, \nu_2). 
$$
Therefore, convergence of measures in Wasserstein distance with rate $e^{-kt}$ (obtained in the main theorems) implies convergence with test function $f$ with rate $e^{-kpt}$. Convergence in $V$-norm for the function $V(x) = e^{\lambda x}$ implies convergence of all moments. Indeed, a simple calculus exercise shows (annd we will make us of it in the proofs) that for every $p \ge 1$ there exists an $a(p) > 0$ such that  $1 + x^p \le a(p)V(x),\, x \ge 0$. Thus for the $f$ as above, we get:
$$
\left|(\nu_1, f) - (\nu_2, f)\right| \le a(p)\norm{\nu_1 - \nu_2}_V.
$$
Therefore, convergence in the norm $\norm{\cdot}_V$ with rate $e^{-\varkappa t}$ implies convergence with test function $f$ with the same rate. To find whether convergence in the $\norm{\cdot}_V$-norm or Wasserstein distance is faster, we need to compare $kp$ with $\varkappa$. Convergence rate in Theorem~\ref{thm:1} is with the same rate as in \cite{Sar16} for the norm  $\norm{\cdot}_V$. However, in Theorem~\ref{thm:2} convergence rate can be faster than in \cite{Sar16} for the norm 
$\norm{\cdot}_V$, if only $K > k(\lambda)$.

\begin{example} For $G = 0$ (including Example 1), Theorem~\ref{thm:2} does not improve upon Theorem~\ref{thm:1}. 
But if $G < 0$, then Theorem~\ref{thm:2} can give a better rate of convergence than Theorem~\ref{thm:1}. For example, for reflected Ornstein-Unlenbeck process (without jumps)
$$
\mathrm{d}X(t) = -a(m+ X(t))\,\mathrm{d}t + \sigma\,\mathrm{d}W(t) + \mathrm{d}\ell(t),\, t \ge 0,
$$
with constants $a, m, \sigma > 0$. Then $G := -a$, and we can take 
$$
\lambda = \frac{am}{\sigma^2},\quad k(\lambda) = \frac{a^2m^2}{2\sigma^2},\quad K := \frac{a^2m^2}{2\sigma^2p} + a. 
$$
This convergence rate $pK$ in Wasserstein distance is faster than the rate $k(\lambda)$ of convergence in the norm $\norm{\cdot}_V$ from \cite{Sar16}:
$$
\frac{a^2m^2}{2\sigma^2} + ap = pK > \frac{a^2m^2}{2\sigma^2}.
$$
\end{example}

\section{Proofs}

\subsection{Proof of Theorem~\ref{thm:1}} (a) The idea is the same as in \cite{Walsh, LMT96, Sar16}: Couple two copies $X_1$ and $X_2$ of this processes starting from $x_1$ and $x_2$. Without loss of generality, assume $0 \le x_1 \le x_2$.  It was shown in \cite[Lemma 4.2]{Sar16} that this process is stochastically ordered if the family of measures $(\nu_x)_{x \ge 0}$ is stochastically ordered, which is assumed by Assumption 2. Thus we can assume $X_1(t) \le X_2(t)$ for all $t \ge 0$. Take $\tau := \inf\{t \ge 0\mid X_2(t) = 0\}$, then $X_1(\tau) = 0$ and we can assume $X_1(t) = X_2(t)$ for $t > \tau$. Using calculus, we  compute 
\begin{equation}
\label{eq:calculus}
x^p \le ae^{\lambda x} = aV(x),\quad x \ge 0,\quad a := (pe/\lambda)^p.
\end{equation}
Since $X_2(t) \ge |X_2(t) - X_1(t)|$ for all $t \ge 0$, $X_1(t) = X_2(t)$ for $t \ge \tau$, by~\eqref{eq:calculus}, 
\begin{align}
\label{eq:new}
\begin{split}
\mathbb E&\left|X_1(t) - X_2(t)\right|^p = \mathbb E\left[1_{\{\tau > t\}}\left|X_1(t) - X_2(t)\right|^p\right] \\ & \le \mathbb E\left[X_2^p(t)1_{\{\tau > t\}}\right] \le a\cdot\mathbb E\left[V(X_2(t))1_{\{\tau > t\}}\right],\quad t \ge 0.
\end{split}
\end{align}
From \cite[Section 5, (5.9)]{Sar16}, the following process is an 
$(\mathcal F_t)_{t \ge 0}$-supermartingale:
\begin{equation}
\label{eq:supermart}
\left(e^{k(\lambda)(t\wedge\tau)}V(X_2(t\wedge\tau)),\, t \ge 0\right).
\end{equation}
Therefore, for all $t \ge 0$ we get:
\begin{equation}
\label{eq:old}
e^{k(\lambda)t}\mathbb E\left[V(X_2(t))1_{\{\tau > t\}}\right] \le \mathbb E\left[e^{k(\lambda)(t\wedge\tau)}V(X_2(t\wedge\tau))\right] \le V(x_2).
\end{equation}
Combining~\eqref{eq:new} and~\eqref{eq:old}, we get:
\begin{equation}
\label{eq:W}
\mathbb E\left|X_1(t) - X_2(t)\right|^p \le ae^{-k(\lambda)t}V(x_2) \le ae^{-k(\lambda)t}V(x_1\vee x_2).
\end{equation}
Raising it to the power $1/p$, we get:
\begin{align*}
\mathcal W_p&(P^t(x_1, \cdot), P^t(x_2, \cdot)) \le \left[\mathbb E\left|X_1(t) - X_2(t)\right|^p\right]^{1/p} \\ & \le a^{1/p}\exp\left[-p^{-1}k(\lambda)t\right]V^{1/p}(x_2) = \lambda^{-1}pe\exp\left[p^{-1}\lambda x_2 - p^{-1}k(\lambda)t\right] \\ & \le \lambda^{-1}pe\exp\left[p^{-1}\lambda (x_1\vee x_2) - p^{-1}k(\lambda)t\right].
\end{align*}

(b) We modify~\eqref{eq:new} and again use~\eqref{eq:calculus}:
\begin{align}
\label{eq:new*}
\begin{split}
\mathbb E&\sup\limits_{t \le s \le T}\left|X_1(s) - X_2(s)\right|^{p} = \mathbb E\left[\sup\limits_{t \le s \le T}\left|X_1(s) - X_2(s)\right|^{p}\cdot 1_{\{\tau > t\}}\right] \\ & \le \mathbb E\left[\sup\limits_{t \le s \le T}X_2^{p}(s)\cdot 1_{\{\tau > t\}}\right] \le a\cdot\mathbb E\left[\sup\limits_{t \le s \le T}V(X_2(t))\cdot 1_{\{\tau > t\}}\right].
\end{split}
\end{align}
We modify~\eqref{eq:old} by applying Ville's maximal inequality from \cite[Exercise 4.8.2]{Durrett} to~\eqref{eq:supermart}:
\begin{equation}
\label{eq:old*}
e^{k(\lambda)t}\cdot\mathbb E\left[\sup\limits_{t \le s \le T}V(X_2(t))\cdot 1_{\{\tau > t\}}\right] \le \mathbb E\left[e^{k(\lambda)(t\wedge\tau)}V(X_2(t\wedge\tau))\right] \le V(x_2).
\end{equation}
Apply~\eqref{eq:Skorohod} to combined~\eqref{eq:new*} and~\eqref{eq:old*} and complete the proof.

(c) In~\eqref{eq:W} from the proof of (a), replace $x_1\vee x_2$ with $x_1 + x_2$, and integrate the right-hand side with respect to $x_1 \sim \rho_1$ and $x_2 \sim \rho_2$. Thus $V(x_1\vee x_2) \le V(x_1) + V(x_2)$. Then we get:
\begin{align*}
\int_0^{\infty}&\int_0^{\infty}[V(x_1)+V(x_2)]\,\rho_1(\mathrm{d}x_1)\,\rho_2(\mathrm{d}x_2)  = \left(\rho_1, V\right) + \left(\rho_2, V\right). 
\end{align*}
Then raise this estimate to the power $1/p$, and complete the proof.

(d) Apply (c) to $\rho_1 = \pi$ and $\rho_2 = \delta_x$, with $\delta_xP^t = P^t(x, \cdot)$ and $\pi P^t = \pi$ for $t \ge 0$.

\subsection{Proof of Theorem~\ref{thm:2}} Without loss of generality, we can assume $x_1 \le x_2$. The following lemma is proved at the end of this subsection.

\begin{lemma} We can couple $X_1$ and $X_2$ starting from $X_1(0) = x_1$ and $X_2(0) = x_2$ so that
$$
0 \le X_2(t) - X_1(t) \le (x_2 - x_1)e^{Gt},\quad t < \tau := \inf\{t \ge 0\mid X_2(t) = 0\}.
$$
\label{lemma:temp}
\end{lemma}

(a) Adapt the proof of Theorem~\ref{thm:1}. Use~\eqref{eq:new} and Lemma~\ref{lemma:temp}:
$$
\mathcal W^p_p(P^t(x_1, \cdot), P^t(x_2, \cdot)) \le \mathbb E\left[\left(X_2(t) - X_1(t)\right)^p1_{\{\tau > t\}}\right]  \le \mathbb E\left[(x_2 - x_1)^pe^{Gpt}1_{\{\tau > t\}}\right].
$$
Since $V(x) \ge 1$ for $x \ge 0$, we get:
$$
e^{Gpt}\cdot\mathbb E\left[1_{\{\tau > t\}}\right] \le e^{Gpt}\cdot \mathbb E\left[1_{\{\tau > t\}}V(X_2(t))\right] \le 
e^{(Gp - k(\lambda))t}\cdot \mathbb E\left[e^{k(\lambda)(t\wedge\tau)}V(X_2(t\wedge\tau))\right].
$$
By the supermartingale property of $(e^{k(\lambda)(t\wedge\tau)}V(t\wedge\tau),\, t \ge 0)$, we get:
$$
\mathbb E\left[e^{k(\lambda)(t\wedge\tau)}V(X_2(t\wedge\tau))\right] \le V(x_2). 
$$
Recall the definition of $K$. Raising  this inequality to the power $1/p$, we complete the proof. 

(b) The proof is similar to the one from (a). Using~\eqref{eq:Skorohod}, we get: 
$$
\mathcal W^p_p(\mathbf{P}^{[t, T]}(x_1, \cdot), \mathbf{P}^{[t, T]}(x_2, \cdot)) \le \mathbb E\,\sup\limits_{t \le s \le T}|X_1(t) - X_2(t)|^p.
$$
Since $G \le 0$, we get from Lemma~\ref{lemma:temp}:
\begin{equation}
\label{eq:12}
\sup\limits_{t \le s \le T}|X_1(t) - X_2(t)|^p \le \sup\limits_{t \le s \le T}e^{Gps}|x_1 - x_2|^p = e^{Gpt}|x_1 - x_2|^p.
\end{equation}
Taking expectation in~\eqref{eq:12}, applying~\eqref{eq:Skorohod}, raising to the power $1/p$, we complete the proof. 

(c, d) The proofs are similar to that of (c, d) from Theorem~\ref{thm:1}. 

\smallskip

{\it Proof of Lemma~\ref{lemma:temp}.} We use {\it piecing out}, a classic method for adding jumps to continuous processes, described in \cite{Piece, Sar16}. Since the intensity $\Lambda := \nu_x(\mathbb R_+)$ is constant, we can couple the processes with simultaneous jumps $0 = \tau_0 < \tau_1 < \tau_2 < \ldots$ as Poisson process with rate $\Lambda$. Let us couple them so that the following is true:
\begin{equation}
\label{eq:ordering}
X_1(t) \le X_2(t),\quad t \ge 0,
\end{equation}
and for $k = 0, 1, \ldots$, if $\tau > \tau_{k+1}$, then
\begin{align}
\label{eq:diff}
\begin{split}
\mathrm{d}&(X_2(t) - X_1(t)) \le G(X_2(t) - X_1(t))\,\mathrm{d}t,\, \tau_k \le t < \tau_{k+1};\\
&X_2(\tau_{k+1}) - X_2(\tau_{k+1}-) \le X_1(\tau_{k+1}) - X_1(\tau_{k+1}-).
\end{split}
\end{align}
Assume we proved~\eqref{eq:ordering} and ~\eqref{eq:diff}. Use induction over $k$ to show the statement of the lemma for $\tau_k < t \le \tau_{k+1}$, assuming this for $t = \tau_k$. By Gronwall's lemma, since $X_2 - X_1$ is continuous,
$$
X_2(t) - X_1(t) \le e^{G(t - \tau_k)}(X_2(\tau_k) - X_1(\tau_k)),\, \tau_k \le t < \tau_{k+1}.
$$
Combining this observation with the second result in~\eqref{eq:diff}, we complete the proof of induction step and with it the proof of Lemma~\ref{lemma:temp}. To prove~\eqref{eq:ordering} and the first line in~\eqref{eq:diff}, for $t \in (\tau_k, \tau_{k+1}]$, use induction over $k = 0, 1, \ldots$ During $(\tau_k, \tau_{k+1})$, $X_1$ and $X_2$ behave as copies of a reflected diffusion with drift $g$ and diffusion $\sigma^2$. Thus
$$
\mathrm{d}X_i(t) = g(X_i(t))\,\mathrm{d}t + \sigma\,\mathrm{d}W(t) + \mathrm{d}\ell_i(t),\, i = 1, 2,
$$
where $\ell_i$ is a continuous nondecreasing process which can increase only when $X_i(t) = 0$, and $\ell_i(0) = 0$. We coupled them with the same driving Brownian motion $W$. Since $X_1(\tau_k) \le X_2(\tau_k)$ by the induction hypothesis, we get~\eqref{eq:ordering} for $t \in (\tau_k, \tau_{k+1})$. For $t < \tau$, $\ell_2(t) = 0$, and $\ell_1$ is nondecreasing, thus $\mathrm{d}\ell_1(t) \ge 0$. Assumption 4 about the function $g$, together with $X_1(t) \le X_2(t)$, implies $g(X_2(t)) - g(X_1(t)) \le G(X_2(t) - X_1(t))$. Therefore,
\begin{align*}
\mathrm{d}&(X_2(t) - X_1(t)) = \left[g(X_2(t)) - g(X_1(t))\right]\,\mathrm{d}t - \mathrm{d}\ell_1(t) \\ &  \le \left[g(X_2(t)) - g(X_1(t))\right]\,\mathrm{d}t \le G\left[X_2(t) - X_1(t)\right]\,\mathrm{d}t.
\end{align*}
This proves the first line in~\eqref{eq:diff}. Letting $t \uparrow \tau_{k+1}$ in $X_1(t) \le X_2(t)$, we get: $y_1 := X_1(\tau_{k+1}-) \le y_2 := X_2(\tau_{k+1}-)$, By Assumption 2 and Assumption 4 for jump measures, there is a coupling $(z_1, z_2)$ of normalized jump measures:
$$
z_1 \sim \Lambda^{-1}\nu_{y_1},\, z_2 \sim \Lambda^{-1}\nu_{y_2},\quad z_2 \le z_1,\, z_1 +y_1 \le z_2 + y_2.
$$
Let $X_i(\tau_{k+1}) := y_i+z_i$, $i = 1, 2$ to prove~\eqref{eq:ordering} and the second line in~\eqref{eq:diff} for $t = \tau_{k+1}$.


\begin{thebibliography}{99}

\bibitem{BCG08} \textsc{Dominique Bakry, Patrick Cattiaux, Arnaud Guillin} (2008). Rate of Convergence of Ergodic Continuous Markov Chains: Lyapunov vs Poincare. \textit{Journal of Functional Analysis} \textbf{254}, 727--754. 

\bibitem{Q} \textsc{Yana Belopolskaya, Guodong Pang, Andrey Sarantsev, Yurii Suhov} (2019). Stationary Distributions and Convergence for M/M/1 Queues in Interactive Random Environment. To appear in \textit{Queueing Systems: Theory and Applications.}


\bibitem{Oleg} \textsc{Oleg Butkowski} (2014). Subgeometric Rates of Convergence of Markov Processes in the Wasserstein Metric. \textit{Annals of Applied Probability} \textbf{24}, 526--552.

\bibitem{CG} \textsc{Patrick Cattiaux, Arnaud Guillin} (2014). Semi Log-Concave Markov Diffusions. \textit{Lecture Notes in Mathematics} \textbf{2123}, 231--292, Springer.



\bibitem{DMT95} \textsc{Douglas Down, Sean P. Meyn, Richard L. Tweedie} (1995). Exponential and Uniform Ergodicity of Markov Processes. \textit{Annals of Probability} \textbf{23}, 1671--1691. 

\bibitem{Durrett} \textsc{Richard Durrett} (2019). \textit{Probability: Theory and Examples.} 5th ed., Cambridge University Press.



\bibitem{Risk} \textsc{Pierre-Olivier Goffard, Andrey Sarantsev} (2019). Exponential Convergence Rate of Ruin Probabilities for Level-Dependent Levy-Driven Risk Processes. \textit{Journal of Applied Probability} \textbf{56}, 1244--1268.

\bibitem{Hairer} \textsc{Martin Hairer, Jonathan Mattingly} (2011). Yet Another Look at Harris' Ergodic Theorem for Markov Chains. \textit{Seminar on Stochastic Analysis, Random Fields and Applications} \textbf{6}, 109--117.

\bibitem{Walsh} \textsc{Tomoyuki Ichiba, Andrey Sarantsev} (2019). Convergence and Stationary Distributions for Walsh Diffusions. \textit{Bernoulli} \textbf{25}, 2439--2478. 

\bibitem{Piece} \textsc{Nobuyuki Ikeda, Masao Nayasawa, Shinzo Watanabe} (1966). A Construction of Markov Processes by Piecing Out. \textit{Proceedings of the Japan Academy of Sciences} \textbf{42}, 370--375.

\bibitem{Lindvall83} \textsc{Torgny Lindvall} (1983). On Coupling of Diffusion Processes. \textit{Journal of Applied Probability} \textbf{20}, 82--93.

\bibitem{LindvallBook} \textsc{Torgny Lindvall} (2002). \textit{Lectures on the Coupling Method.} Dover. 

\bibitem{Lindvall86} \textsc{Torgny Lindvall, L. C. G. Rogers} (1986). Coupling of Multidimensional Diffusions by Reflection. \textit{Annals of Probability} \textbf{14}, 860--872.

\bibitem{LMT96} \textsc{Robert B. Lund, Sean P. Meyn, Richard L. Tweedie} (1996). Computable Exponential Convergence Rates for Stochastically Ordered Markov Processes. \textit{Annals of Applied Probability} \textbf{6}, 218--237. 

\bibitem{MT93a} \textsc{Sean P. Meyn, Richard L. Tweedie} (1993). Stability of Markovian Processes II. Continuous-Time Processes and Sampled Chains. \textit{Advances in Applied Probability} \textbf{25}, 487--517. 

\bibitem{MT93b} \textsc{Sean P. Meyn, Richard L. Tweedie} (1993). Stability of Markovian Processes III. Foster-Lyapunov Criteria for Continuous-Time Markov Processes. \textit{Advances in Applied Probability} \textbf{25}, 518--548. 

\bibitem{MT94} \textsc{Sean P. Meyn, Richard L. Tweedie} (1994). Computable Bounds for Geometric Convergence Rates of Markov Chains. \textit{Annals of Applied Probability} \textbf{4}, 981--1011. 

\bibitem{Pal} \textsc{Soumik Pal} (2012). Concentration for Multidimensional Diffusions and their Boundary Local
Times. \textit{Probability Theory and Related Fields} \textbf{154}, 225--254.

\bibitem{PS19} \textsc{Soumik Pal and Andrey Sarantsev} (2019). A Note on Transportation Cost Inequalities for Diffusions with
Reflections. \textit{Electronic Communications in Probability} \textbf{24}, \# 21.

\bibitem{MCMC} \textsc{Jeffrey S. Rosenthal} (1995). Minorization Conditions and Convergence Rates for Markov Chain Monte Carlo. \textit{Journal of the American Statistical Association} \textbf{90}, 558--566. 

\bibitem{Samson} \textsc{Paul-Marie Samson} (2000). Concentration of Measure Inequalities for Markov Chains and
$\Phi$-Mixing Processes. \textit{Annals of Probability} \textbf{28}, 416--461.

\bibitem{Sar16} \textsc{Andrey Sarantsev} (2016). Explicit Rates of Exponential Convergence for Reflected Jump-Diffusions on the Half-Line. \textit{Latin American Journal of Probability \& Mathematical Statistics} \textbf{13}, 1069--1093. 

\bibitem{Villani} \textsc{Cedric Villani} (2008). \textit{Optimal Transport: Old and New.} \textit{Grundlehren der Mathematischen Wissenschaften} \textbf{338}, Springer. 

\end{thebibliography}
\end{document}